\documentstyle[12pt]{article}
\begin {document}

\makeatletter
\@addtoreset{equation}{section}

\makeatother
\makeatletter

\begin{center}{\Large \bf Introduction to vertex operator
 algebras II}\end{center}

\vspace{0.5cm}
\begin{center}{\bf Hai-sheng Li\\
 Department of Mathematics\\
University of California\\
Santa Cruz, CA 95064}\end{center}
\vspace{0.5cm}

\baselineskip=24pt

\section{Introduction}
This is the second of three lectures on introduction to vertex operator
algebras. In this lecture, we shall continue Professor Dong's lecture to
present more fundamental properties of vertex operator algebras.

{}From the mathematical point of view,  a vertex operator algebra formally
resembles a Lie algebra because the Jacobi identity is used as
 one of the main axioms.
For the Lie algebra aspect of vertex operator algebras, the notion of
contragredient module [FHL] and the notion of tensor product ([HL1-4], [Li4])
have been developed. On the other hand, from the physical point of view, a
vertex operator algebra
looks like a commutative associative algebra with identity because
roughly speaking, a vertex operator algebra is a sort of quantization of
the commutative associative algebra of observables in conformal field theory
 [BPZ]. For the associative
algebra  aspect of vertex operator algebras, it has been proved [FHL] that
the tensor product of any finitely many vertex operator algebras has a natural
vertex operator algebra structure. In concrete examples, for a fixed level
$\ell$, one of the generalized Verma modules, called the vacuum representation
 for any affine Lie algebra $\tilde{{\bf g}}$, has a natural
vertex operator algebra structure ([FZ], [Li2], [Lia]) and the universal
enveloping  algebra [FZ] of
the vertex operator algebra is a certain completion of the universal
enveloping algebra of $\tilde{{\bf g}}$. This fact together with some
facts on tensor products ([HL1-4], [Li4], [KL0-2],...) strongly indicates
that a vertex operator algebra is analogous to a quasi-Hopf algebra or a
quantum group.

In this lecture, we shall discuss certain analogies between vertex operator
algebras and classical algebras such as commutative associative algebras and
Lie algebras. More specifically, we shall discuss commutativity and
associativity and we give an analogue of the endomorphism ring for vertex
operator algebras.
All the materials presented here are taken from [DL], [FLM], [FHL] and [Li2].

{\large \bf Acknowledgment} We'd like to thank Professor Miyamoto for his
hospitality during
the short visit in Japan.

\section{Commutativity and associativity}
In this section, we shall present two versions of commutativity and
associativity, which are given in terms of formal variables
and in terms of analytic functions, respectively. In [FHL] and [FLM], it was
proved that
the Jacobi identity is equivalent to the commutativity and associativity in
terms of
``matrix-coefficient'' and formal variables. This technique was further
employed in [DL]
for a more general
structure, called generalized vertex operator algebras, or even more general
structure, called
abelian intertwining algebras, where
the version of the commutativity and the associativity in terms of formal
variables
without using ``matrix-coefficient'' was first noticed.
Since then  it has been realized that the matrix-coefficient-free formal
variable technique
is very useful some time. For instance, when dealing with vertex algebras
(recalled later) we don't
have a grading available to define the ``matrix-coefficient,'' this
matrix-coefficient-free version
is very useful. (See [A], [Gu] for more applications.)

\subsection{Commutativity and associativity in terms of formal variables}
A vertex operator algebra [FLM] has been already recalled in the first lecture.
Now we briefly recall the definition of a vertex algebra from [B]. A vertex
algebra satisfies all
the  axioms for a vertex operator algebra except that we don't assume the
existence of
the Virasoro algebra so that the ${\bf Z}$-grading is not assumed.

{\bf Proposition 2.1. } {\it Let $V$ be a vertex (operator) algebra and
let $a$ and $b$ be two elements of $V$. Then there is a nonnegative integer
$k$ such that}
\begin{eqnarray}
(z_{1}-z_{2})^{k}Y(a,z_{1})Y(b,z_{2})=(z_{1}-z_{2})^{k}Y(b,z_{2})Y(a,z_{1}).
\end{eqnarray}

{\bf Proof.} In order to obtain (2.1), we need to force the term on the right
hand side to vanish. For any nonnegative integer $m$, taking
${\rm Res}_{z_{0}}z_{0}^{m}$ of the Jacobi identity, we obtain:
\begin{eqnarray}
& &(z_{1}-z_{2})^{m}[Y(a,z_{1}),Y(b,z_{2})]\nonumber\\
&=&(z_{1}-z_{2})^{m}\left(Y(a,z_{1})Y(b,z_{2})-
Y(b,z_{2})Y(a,z_{1})\right)\nonumber\\
&=&{\rm Res}_{z_{0}}z_{2}^{-1}\delta\left(
\frac{z_{1}-z_{0}}{z_{2}}\right)z_{0}^{m}Y(Y(a,z_{0})b,z_{2})\nonumber \\
&=&\sum_{i=0}^{\infty}\frac{1}{i!}\left(\left({\partial\over\partial z_{2}}
\right)^{i}z_{1}^{-1}\delta\left(\frac{z_{2}}{z_{1}}\right)\right)
Y(a_{m+i}b,z_{2}).
\end{eqnarray}
Let $k$ be a nonnegative integer such that
$a_{n}b=0$ for $n\ge k$. Then setting $m=k$ we obtain the commutativity (2.1).
$\;\;\;\;\Box$

{\bf Proposition 2.2. } {\it Let $V$ be a vertex (operator) algebra and
let $a$ and $c$ be two elements of $V$. Then there is a nonnegative integer
$k$ such that for any } $b\in V$ {\it we have}
\begin{eqnarray}
(z_{0}+z_{2})^{k}Y(Y(a,z_{0})b,z_{2})c=(z_{0}+z_{2})^{k}Y(a,z_{0}+z_{2})
Y(b,z_{2})c.
\end{eqnarray}

{\bf Proof.} Similar to the case for commutativity, in order to obtain the
associativity we force the second term on the left hand side to be zero.
Taking ${\rm Res}_{z_{1}}$ of the Jacobi identity, we obtain the following
{\it iterate formula}:
\begin{eqnarray}
& &Y(Y(a,z_{0})b,z_{2})\nonumber \\
&=&{\rm Res}_{z_{1}}\left(z_{0}^{-1}\delta\left(\frac{z_{1}-z_{2}}{z_{0}}
\right)Y(a,z_{1})Y(b,z_{2})-z_{0}^{-1}\delta\left(\frac{z_{2}-z_{1}}{-z_{0}}
\right)Y(b,z_{2})Y(a,z_{1})\right)\nonumber \\
&=&Y(a,z_{0}+z_{2})Y(b,z_{2})-Y(b,z_{2})(Y(a,z_{0}+z_{2})-Y(a,z_{2}+z_{0})).
\end{eqnarray}
For any $c\in V$, let $m$ be a positive integer such that
$z^{m}Y(a,z)c$ involves only positive powers of $z$, so that
\begin{eqnarray}
(z_{0}+z_{2})^{m}(Y(a,z_{0}+z_{2})-Y(a,z_{2}+z_{0}))c=0.
\end{eqnarray}
Then we obtain the associativity (2.3).$\;\;\;\;\Box$

The following theorem has been proved in [FHL] and [FLM], and in [DL] for a
more general
context. Here,  as in [Li2] we give a slightly different proof without using
 ``matrix-coefficient.''

{\bf Theorem 2.3.} {\it The Jacobi identity is equivalent to the
commutativity together with the associativity.}

{\bf Proof.} We only need to prove that the Jacobi identity follows from
commutativity and associativity. Choosing a nonnegative integer $k$ such that
$a_{m}c=0$ for all $m\ge k$, we get
\begin{eqnarray}
& &z_{0}^{k}z_{1}^{k}\left(z_{0}^{-1}\delta\left(\frac{z_{1}-z_{2}}{z_{0}}
\right)Y(a,z_{1})Y(b,z_{2})c-z_{0}^{-1}\delta\left(
\frac{z_{2}-z_{1}}
{-z_{0}}\right)Y(b,z_{2})Y(a,z_{1})c\right)\nonumber\\
&=&z_{0}^{-1}\delta\left(\frac{z_{1}-z_{2}}{z_{0}}\right)z_{1}^{k}
(z_{1}-z_{2})^{k}Y(a,z_{1})Y(b,z_{2})c\nonumber\\
& &-z_{0}^{-1}\delta\left(\frac{z_{2}-z_{1}}{-z_{0}}\right)
z_{1}^{k}(z_{1}-z_{2})^{k}Y(b,z_{2})Y(a,z_{1})c\nonumber\\
&=&z_{2}^{-1}\delta\left(\frac{z_{1}-z_{0}}{z_{2}}\right)
\left(z_{1}^{k}(z_{1}-z_{2})^{k}Y(b,z_{2})Y(a,z_{1})c\right)\nonumber\\
&=&z_{2}^{-1}\delta\left(\frac{z_{1}-z_{0}}{z_{2}}\right)
\left(z_{0}^{k}(z_{0}+z_{2})^{k}Y(b,z_{2})Y(a,z_{2}+z_{0})c\right).
\end{eqnarray}
Since $a_{m}c=0$ for all $m\ge k$,
$(z_{0}+z_{2})^{k}Y(a,z_{2}+z_{0})c$ involves only nonnegative powers
of $(z_{2}+z_{0})$, so that
\begin{eqnarray}
z_{0}^{k}(z_{0}+z_{2})^{k}Y(b,z_{2})Y(a,z_{2}+z_{0})c
=z_{0}^{k}(z_{0}+z_{2})^{k}Y(b,z_{2})Y(a,z_{0}+z_{2})c.
\end{eqnarray}
Therefore
\begin{eqnarray}
& &z_{0}^{k}z_{1}^{k}\left(z_{0}^{-1}\delta\left(\frac{z_{1}-z_{2}}{z_{0}}
\right)Y(a,z_{1})Y(b,z_{2})c-z_{0}^{-1}\delta\left(
\frac{z_{2}-z_{1}}
{-z_{0}}\right)Y(b,z_{2})Y(a,z_{1})c\right)\nonumber\\
&=&z_{2}^{-1}\delta\left(\frac{z_{1}-z_{0}}{z_{2}}\right)
\left(z_{0}^{k}(z_{0}+z_{2})^{k}Y(b,z_{2})Y(a,z_{0}+z_{2})c\right)
\nonumber\\
&=&z_{2}^{-1}\delta\left(\frac{z_{1}-z_{0}}{z_{2}}\right)\left(z_{0}^{k}
(z_{0}+z_{2})^{k}Y(a,z_{0}+z_{2})Y(b,z_{2})c\right)\nonumber\\
&=&z_{2}^{-1}\delta\left(\frac{z_{1}-z_{0}}{z_{2}}\right)\left(z_{0}^{k}
(z_{0}+z_{2})^{k}Y(Y(a,z_{0})b,z_{2})c\right)\nonumber\\
&=&z_{2}^{-1}\delta\left(\frac{z_{1}-z_{0}}{z_{2}}\right)\left(z_{0}^{k}
z_{1}^{k}Y(Y(a,z_{0})b,z_{2})c\right)\nonumber\\
&=&z_{0}^{k}z_{1}^{k}z_{2}^{-1}\delta\left(\frac{z_{1}-z_{0}}{z_{2}}\right)
Y(Y(a,z_{0})b,z_{2})c.
\end{eqnarray}
Then the Jacobi identity follows. $\;\;\;\;\Box$

{\bf Remark 2.4}. Notice that the commutativity (2.1) is
really a commutativity for ``left multiplications.'' The exact analogue of the
classical commutativity of product is the {\it skew-symmetry} ([B], [FHL]):
$Y(a,z)b=e^{zL(-1)}Y(b,-z)a$ for any $a,b\in V$. Let
$A$ be a classical  algebra with a right identity 1 and denote by $\ell _{a}$
 the left multiplication by an element $a$. Suppose that
$\ell_{a}\ell_{b}=\ell_{b}\ell_{a}$ for any $a,b\in A$. Then
 $a(bc)=b(ac)$ for any $a,b,c\in A$. Setting
$c=1$, we obtain the commutativity $ab=ba$. Furthermore, we obtain the
associativity:
\begin{eqnarray}
a(cb)=a(bc)=b(ac)=(ac)b\;\;\;\mbox{ for any }a,b,c\in A.
\end{eqnarray}
Therefore, $A$ is a commutative associative algebra. This classical fact
suggests that the commutativity (2.1) together with the vacuum property
implies the associativity (2.3).
Therefore, the commutativity implies the Jacobi identity. This has been
proved in slightly different contexts in many references such as [FLM], [FHL],
[Go], [DL]
and [Li2].

{\bf Theorem 2.5}. {\it In the definition of vertex operator algebra, the
Jacobi identity can be equivalently replaced by the commutativity.}

{\bf Proof.} The proof, which consists of three steps,
is exactly an analogue of the argument given in Remark 2.4. This is essentially
the same
proof in [FHL] more general [DL].  However,
except that we use matrix-coefficient-free technique we
will prove the skew-symmetry first so that we will not
use the product of three vertex operators (see [FHL, [DL]).

(1) The skew-symmetry holds. Let $k$ be a positive
integer such that
$b_{m}a=0$ for all $m\ge k$ and that the commutativity (2.1)
holds. Then
\begin{eqnarray}
& &(z_{1}-z_{2})^{k}Y(a,z_{1})Y(b,z_{2}){\bf 1}\nonumber\\
&=&(z_{1}-z_{2})^{k}Y(b,z_{2})Y(a,z_{1}){\bf 1}\nonumber\\
&=&(z_{1}-z_{2})^{k}Y(b,z_{2})e^{z_{1}L(-1)}a\nonumber \\
&=&(z_{1}-z_{2})^{k}e^{z_{1}L(-1)}Y(b,z_{2}-z_{1})a.
\end{eqnarray}
Since $(z_{1}-z_{2})^{k}Y(b,z_{2}-z_{1})a$ involves only nonnegative
powers of $(z_{2}-z_{1})$, we may set
$z_{2}=0$. Thus
\begin{eqnarray}
z_{1}^{k}Y(a,z_{1})b=z_{1}^{k}e^{z_{1}L(-1)}Y(b,-z_{1})a.
\end{eqnarray}
Multiplying both sides of (2.11) by $z_{1}^{-k}$ we obtain
 $Y(a,z_{1})b=e^{z_{1}L(-1)}Y(b,-z_{1})a$.

(2) The associativity (2.3) holds.
For any $a,c\in V$, let $k$ be a
positive integer such that
the commutativity (2.1) for $(a,c)$ holds. Then for any
$b\in V$, we have:
\begin{eqnarray}
& &(z_{0}+z_{2})^{k}Y(a,z_{0}+z_{2})Y(b,z_{2})c\nonumber\\
&=&(z_{0}+z_{2})^{k}Y(a,z_{0}+z_{2})e^{z_{2}L(-1)}Y(c,-z_{2})b
\nonumber\\
&=&e^{z_{2}L(-1)}(z_{0}+z_{2})^{k}Y(a,z_{0})
Y(c,-z_{2})b\nonumber\\
&=&e^{z_{2}L(-1)}(z_{0}+z_{2})^{k}Y(c,-z_{2})Y(a,z_{0})b\nonumber\\
&=&(z_{0}+z_{2})^{k}Y(Y(a,z_{0})b,z_{2})c.
\end{eqnarray}
It follows from Theorem 2.3 that the Jacobi identity holds. $\;\;\;\;\Box$

The following Proposition 2.6, Corollary 2.7 and Proposition 2.8 are taken
from [Li2].

{\bf Proposition 2.6.} {\it The Jacobi identity in the definition of vertex
operator algebra $V$ can be equivalently replaced by the skew symmetry
and the associativity (2.3).}

{\bf Proof}.  For any
$a,b,c\in V$, let $k$ be a positive integer such that $z^{k}Y(b,z)c$
involves only positive powers of $z$ and that the following
associativities hold:
\begin{eqnarray}
(z_{0}+z_{2})^{k}Y(a,z_{0}+z_{2})Y(b,z_{2})c
&=&(z_{0}+z_{2})^{k}Y(Y(a,z_{0})b,z_{2})c,\nonumber\\
(-z_{0}+z_{1})^{k}Y(b,-z_{0}+z_{1})Y(a,z_{1})c&=&(-z_{0}+z_{1})^{k}
Y(Y(b,-z_{0})a,z_{1})c.
\end{eqnarray}
Then
\begin{eqnarray}
& &z_{1}^{k}z_{2}^{k}\left(z_{0}^{-1}\delta\left(\frac{z_{1}-z_{2}}
{z_{0}}\right)Y(a,z_{1})Y(b,z_{2})c-z_{0}^{-1}\delta\left(
\frac{-z_{2}+z_{1}}
{z_{0}}\right)Y(b,z_{2})Y(a,z_{1})c\right)\nonumber\\
&=&z_{0}^{-1}\delta\left(\frac{z_{1}-z_{2}}
{z_{0}}\right)\left((z_{0}+z_{2})^{k}z_{2}^{k}Y(Y(a,z_{0})b,z_{2})c\right)
\nonumber\\
& &-z_{0}^{-1}\delta\left(\frac{-z_{2}+z_{1}}
{z_{0}}\right)\left(z_{1}^{k}(-z_{0}+z_{1})^{k}Y(Y(b,-z_{0})a,z_{1})c\right)
\nonumber\\
&=&z_{0}^{-1}\delta\left(\frac{z_{1}-z_{2}}
{z_{0}}\right)\left((z_{0}+z_{2})^{k}z_{2}^{k}Y(Y(a,z_{0})b,z_{2})c\right)
\nonumber\\
& &-z_{0}^{-1}\delta\left(\frac{-z_{2}+z_{1}}
{z_{0}}\right)\left(z_{1}^{k}(-z_{0}+z_{1})^{k}Y(e^{-z_{0}L(-1)}Y(a,z_{0})
b,z_{1})c\right)\nonumber\\
&=&z_{0}^{-1}\delta\left(\frac{z_{1}-z_{2}}
{z_{0}}\right)\left((z_{0}+z_{2})^{k}z_{2}^{k}Y(Y(a,z_{0})b,z_{2})c\right)
\nonumber\\
& &-z_{0}^{-1}\delta\left(\frac{-z_{2}+z_{1}}
{z_{0}}\right)\left(z_{1}^{k}(-z_{0}+z_{1})^{k}Y(Y(a,z_{0})b,z_{1}-z_{0})c
\right).
\end{eqnarray}
Since $z_{2}^{k}(z_{0}+z_{2})^{k}Y(Y(a,z_{0})b,z_{2})c=(z_{0}+z_{2})^{k}
Y(a,z_{0}+z_{2})(z_{2}^{k}Y(b,z_{2})c)$ involves only
positive powers of $z_{2}$, by (2.13) we have:
\begin{eqnarray}
& &z_{0}^{-1}\delta\left(\frac{z_{1}-z_{2}}
{z_{0}}\right)\left((z_{0}+z_{2})^{k}z_{2}^{k}Y(Y(a,z_{0})b,z_{2})c\right)
\nonumber\\
&=&z_{0}^{-1}\delta\left(\frac{z_{1}-z_{2}}
{z_{0}}\right)\left(z_{1}^{k}(z_{1}-z_{0})^{k}Y(Y(a,z_{0})b,z_{1}-z_{0})c
\right).
\end{eqnarray}
Thus
\begin{eqnarray}
& &z_{1}^{k}z_{2}^{k}\left(z_{0}^{-1}\delta\left(\frac{z_{1}-z_{2}}
{z_{0}}\right)Y(a,z_{1})Y(b,z_{2})c-z_{0}^{-1}\delta\left(
\frac{-z_{2}+z_{1}}
{z_{0}}\right)Y(b,z_{2})Y(a,z_{1})c\right)\nonumber\\
&=&z_{0}^{-1}\delta\left(\frac{z_{1}-z_{2}}
{z_{0}}\right)\left(z_{1}^{k}(z_{1}-z_{0})^{k}Y(Y(a,z_{0})b,z_{1}-z_{0})c
\right)\nonumber\\
& &-z_{0}^{-1}\delta\left(\frac{-z_{2}+z_{1}}
{z_{0}}\right)\left(z_{1}^{k}(-z_{0}+z_{1})^{k}Y(Y(a,z_{0})b,z_{1}-z_{0})c
\right)\nonumber\\
&=&z_{2}^{-1}\delta\left(\frac{z_{1}-z_{0}}
{z_{2}}\right)\left(z_{1}^{k}(z_{1}-z_{0})^{k}Y(Y(a,z_{0})b,z_{1}-z_{0})c
\right)\nonumber\\
&=&z_{1}^{k}z_{2}^{k}z_{2}^{-1}\delta\left(\frac{z_{1}-z_{0}}{z_{2}}\right)
Y(Y(a,z_{0})b,z_{2})c.
\end{eqnarray}
Multiplying both sides by $z_{1}^{-k}z_{2}^{-k}$, we obtain the Jacobi
identity.$\;\;\;\;\Box$

In [B], Borcherds first defined the notion of vertex algebra with a
set of axioms consisting of the vacuum property, the skew-symmetry
and the iterate formula (2.4). Without assuming the existence of a Virasoro
element in the notion of vertex algebra, one can define
the operator $D$ is defined by $Da=a_{-2}{\bf 1}$ for $a\in V$.
Therefore, we have:

{\bf Corollary 2.7}. {\it Borcherds' definition [B] and FLM's definition
for a vertex algebra are equivalent.}

{\bf Proposition 2.8}. {\it Let $V$ be a
vertex algebra. Then the Jacobi identity of a $V$-module can
be equivalently replaced by the associativity (2.3).}

\subsection{Commutativity and associativity in terms of analytic functions}
In the last subsection, we consider a vertex operator $Y(a,z)$ as a generating
function of a sequence of operators, where $z$ is a formal variable. In this
subsection
we shall consider a vertex operator $Y(a,z)$ as an operator valued functional,
where $z$ will be considered as a nonzero complex number.

Let $V$ be a vertex operator algebra and define
$V'=\oplus_{n\in {\bf Z}}V_{(n)}^{*}$ to be the restricted
 dual of $V$ [FHL]. For any $a, b\in V$, by the definition of a vertex operator
 algebra, $Y(a,z)b$ involves only finitely many negative powers of $z$. Then
 it follows from the definition of $V'$ that for any $ f\in V'$,
$\langle f, Y(a,z)b\rangle$ is a Laurent polynomial in $z$. Therefore, we
may consider $z$ as a nonzero complex number so that
$\langle f, Y(a,z)b\rangle$ is a rational function of $z$. Next, we consider
the following {\it $n$-point functions}:
\begin{eqnarray}
\langle f,Y(a_{1},z_{1})Y(a_{2},z_{2})\cdots Y(a_{n},z_{2})c\rangle,
\end{eqnarray}
The following
Proposition 2.9 was taken from [FLM] and [FHL] with a slightly different proof.

{\bf Proposition 2.9} {\it Let $V$ be a vertex operator algebra and let
$a,b,c\in V, f\in V'$. Then (I) (Rationality)}
\begin{eqnarray}
\langle f,Y(a,z_{1})Y(b,z_{2})c\rangle,
\end{eqnarray}
{\it as a formal series converges in the domain of}
 $|z_{1}|>|z_{2}|>0$ {\it to a rational
 function
$g(z_{1},z_{2})=h(z_{1},z_{2})z_{1}^{m}z_{2}^{n}(z_{1}-z_{2})^{k}$, where
 $h(z_{1},z_{2})$ is a polynomial in $z_{1}$ and $z_{2}$, and $m,n,k$ are
integers and  $k$ only depends on $a$ and $b$.}

{\it (II) (Commutativity)}
\begin{eqnarray}
\langle f,Y(b,z_{2})Y(a,z_{1})c\rangle,
\end{eqnarray}
{\it as a formal series converges in the domain of}
 $|z_{2}|>|z_{1}|>0$ {\it to the same rational
 function  $g(z_{1},z_{2})$ as that in (I).}

{\it (III) (Associativity)}
\begin{eqnarray}
\langle f,Y(Y(a,z_{0})b,z_{2})c\rangle,
\end{eqnarray}
{\it as a formal series converges in the domain of}
 $|z_{2}+z_{0}|>|z_{2}|>|z_{0}|>0$ {\it to the rational
 function  $g(z_{2}+z_{0},z_{2})$, where $g(z_{1},z_{2})$ is the same as that
 in (I) and (II).}

{\bf Proof.} For any $a,b\in V$, it follows from the first version of
commutativity that there is a nonnegative integer $k$ such that the
commutativity (2.1) holds. It follows from the definition of $V'$ and the
meromorphic condition on vertex operators that the matrix-coefficient
$\langle f,(z_{1}-z_{2})^{k}Y(a,z_{1})Y(b,z_{2})c\rangle$ involves only
finitely many negative powers of $z_{2}$ and finitely many positive powers
of $z_{1}$. Similarly, $\langle f,Y(b,z_{2})Y(a,z_{1})c\rangle$ involves only
finitely many negative powers of $z_{1}$ and
finitely many positive powers of $z_{2}$. Therefore the common formal series
\begin{eqnarray}
\langle f,(z_{1}-z_{2})^{k}Y(a,z_{1})Y(b,z_{2})c\rangle=
\langle f,(z_{1}-z_{2})^{k}Y(b,z_{2})Y(a,z_{1})c\rangle
\end{eqnarray}
involves only finitely many  negative and positive powers of both $z_{1}$ and
 $z_{2}$. Consequently, it gives a rational function in the form of
$z_{1}^{m}z_{2}^{n}h(z_{1},z_{2})$, where $h(z_{1},z_{2})$ is a polynomial in
 $z_{1}$ and $z_{2}$. Then both (I) and (II) have been proved.
Similarly one can prove (III).$\;\;\;\;\Box$

{\bf Proposition [DL]} {\it The two versions of commutativity and
associativity are equivalent.}

\section{An analogue of the endomorphism ring for VOA}
In this section, we shall introduce what we call ``local systems of
vertex operators'' for any vector space $M$ and we prove that any
local system has a natural vertex algebra structure
with $M$ as a module. Furthermore, we prove that for a fixed vertex
algebra $V$, giving a $V$-module $M$ is equivalent to giving a
vertex algebra homomorphism from $V$ to some local system of
vertex operators on $M$. This whole section was taken from [Li2].
An analogue of the homomorphism module for vertex
operator algebras has been also developed in [Li2] and the notion of local
systems of vertex operators and applications have been generalized to the
notion of local systems of twisted vertex operators in [Li3].

{\bf Definition 3.1.} Let $M$ be any
vector space.  A {\it weak vertex operator} on $M$ is a
formal series $a(z)=\sum_{n \in {\bf
Z}}a_{n}z^{-n-1} \in ({\rm End}\:M)[[z,z^{-1}]]$ such that
\begin{eqnarray}
a(z)u\in M((z))\;\;\;\mbox{for any }u\in M.
\end{eqnarray}
That is, $a_{n}u=0$ for $n$ sufficiently large.
Let $(M,d)$ be a pair consisting of a
vector space $M$ and an endomorphism $d$ of $M$.
A {\it weak vertex operator
on $(M,d)$} is a weak vertex operator $a(z)$ on $M$ such that
\begin{eqnarray}
[d,a(z)]=a'(z) \left(={d\over dz}a(z)\right).
\end{eqnarray}

Denote by $F(M)$ (resp. $F(M,d)$) the space of all
weak vertex operators on $M$ (resp. $(M,d)$).

By definition, it is clear that if $a(z)$ is a weak
vertex operator on
$M$ (resp. $(M,d)$), the formal derivative $a'(z)$ is also a weak
vertex operator on $M$ (resp. $(M,d)$). Then we have an endomorphism
$\displaystyle{D={d\over dz}}$ for both $F(M)$ and $F(M,d)$.

{\bf Definition 3.2.} Let $M$ be a restricted $Vir$-module of central
charge $\ell$. A
weak vertex operator $a(z)$ on $(M,L(-1))$ is
said to be of  {\it  weight} $h\in {\bf C}$  if it
satisfies the following condition:
\begin{eqnarray}
[L(0),a(z)]=ha(z)+za'(z).
\end{eqnarray}
Denote by $F(M,L(-1))_{(h)}$ the space of weak vertex operators
on $(M, L(-1))$ of weight $h$ and set
\begin{eqnarray}
F^{o}(M,L(-1))=\oplus_{h\in {\bf C}}F(M,L(-1))_{(h)}.
\end{eqnarray}

{\bf Remark 3.3}. For any vector space $M$, the identity operator
$I(z)={\rm id}_{M}$ is a weak vertex operator on $M$.
Let $M$ be a restricted $Vir$-module. Then $I(z)={\rm
id}_{M}$ is a weak vertex operator on $(M, L(-1))$ of weight zero and
 $L(z)=\sum_ {n \in {\bf Z}}L(n)z^{-n-2}$
is a weak vertex operator on $(M,L(-1))$ of weight two.
If $a(z)$ is a weak vertex
operator on $(M,L(-1))$ of weight $h$, then
$\displaystyle{a'(z)={d\over dz}a(z)}$ is
a weak vertex operator of weight $h +1$.

{\bf Lemma 3.4}. {\it Let $M$ be a vector space and let $a(z)$
and $b(z)$ be weak vertex operators
on $M$. For any integer $n$, set}
\begin{eqnarray}
a(z)_{n}b(z)={\rm Res}_{z_{1}}\left((z_{1}-z)^{n}a(z_{1})b(z)
-(-z+z_{1})^{n}b(z)a(z_{1})\right).
\end{eqnarray}
{\it Then $a(z)_{n}b(z)$ is a weak
vertex operator.}

{\bf Proof.} For any $u\in M$, by definition we have
\begin{eqnarray}
& &(a(z)_{n}b(z))u\nonumber\\
&=&{\rm Res}_{z_{1}}\left((z_{1}-z)^{n}a(z_{1})b(z)u
-(-z+z_{1})^{n}b(z)a(z_{1})u\right)\nonumber\\
&=&\sum_{k=0}^{\infty}\left(\begin{array}{c}n\\k\end{array}\right)\left(
(-z)^{k}a_{n-k}b(z)u-(-z)^{n-k}b(z)a_{k}u\right).
\end{eqnarray}
It is easy to see that $(a(z)_{n}b(z))u\in M((z))$. Therefore,
$a(z)_{n}b(z)$ is a weak vertex operator on $M$.$\;\;\;\;\Box$

{\bf Definition 3.5}. Let $M$ be a vector space and let
$a(z)$ and $b(z)$ be weak vertex operators on $M$. Then we define
\begin{eqnarray}
& &Y(a(z),z_{0})b(z)\nonumber\\
&=&: \sum_{n\in {\bf Z}}a(z)_{n}b(z) z_{0}^{-n-1}\nonumber\\
&=&{\rm Res}_{z_{1}}
\left(z_{0}^{-1}\delta\left(\frac{z_{1}-z}{z_{0}}\right)a(z_{1})
b(z)-z_{0}^{-1}
\delta\left(\frac{z-z_{1}}{-z_{0}}\right)b(z)a(z_{1})\right).
\end{eqnarray}
Extending the definition bilinearly, we obtain a linear map
\begin{eqnarray}
Y(\cdot,z_{0}) :& & F(M)\rightarrow ({\rm End}
F(M))[[z_{0},z_{0}^{-1}]];\nonumber\\
 & &a(z)\mapsto Y(a(z),z_{0}).
\end{eqnarray}

{\bf Lemma 3.6}. {\it For any $a(z)\in F(M)$, we have}
\begin{eqnarray}
& &Y(I(z),z_{0})a(z)=a(z);\\
& &Y(a(z),z_{0})I(z)=e^{z_{0}{\partial\over\partial
z}}a(z)\;(=a(z+z_{0})).
\end{eqnarray}

{\bf Proof}. By definition, we have:
\begin{eqnarray*}
& &Y(I(z),z_{0})a(z)\\
&=&{\rm Res}_{z_{1}}\left(z_{0}^{-1}\delta\left(\frac{z_{1}-z}{z_{0}}
\right)a(z)-z_{0}^{-1}\delta\left(\frac{-z+z_{1}}{z_{0}}\right)
a(z)\right)\\
&=&{\rm Res}_{z_{1}}z^{-1}\delta\left(\frac{z_{1}-z_{0}}{z}\right)
a(z)\\
&=&a(z)
\end{eqnarray*}
and
\begin{eqnarray}
& &Y(a(z),z_{0})I(z)\nonumber\\
&=&{\rm Res}_{z_{1}}
\left(z_{0}^{-1}\delta\left(\frac{z_{1}-z}{z_{0}}\right)a(z_{1})I(z)-
z_{0}^{-1}\delta\left(\frac{-z+z_{1}}{z_{0}}\right)I(z)a(z_{1})
\right)\nonumber\\
&=&{\rm Res}_{z_{1}}
z^{-1}\delta\left(\frac{z_{1}-z_{0}}{z}\right)a(z_{1})\nonumber\\
&=&{\rm Res}_{z_{1}}z_{1}^{-1}\delta\left(\frac{z+z_{0}}
{z_{1}}\right)a(z_{1})\nonumber\\
&=&{\rm Res}_{z_{1}}z_{1}^{-1}\delta\left(\frac{z+z_{0}}
{z_{1}}\right)a(z+z_{0})\nonumber\\
&=&a(z+z_{0})\nonumber\\
&=&e^{z_{0}{\partial\over\partial z}}a(z).\;\;\;\;\Box
\end{eqnarray}

{\bf Lemma 3.7}. {\it Let $M\in ob C$ and  $a(z),b(z)\in F(M)$.
Then we have}
\begin{eqnarray}
{\partial\over \partial z_{0}}Y(a(z),z_{0})b(z)=
Y(D(a(z)),z_{0})b(z)=[D,Y(a(z),z_{0})]b(z).
\end{eqnarray}

{\bf Proof.} By definition, we have
\begin{eqnarray}
& &{\partial\over \partial z_{0}}Y(a(z),z_{0})b(z)\nonumber\\
&=&{\rm Res}_{z_{1}}\left({\partial \over \partial z_{0}}\left( z_{0}^{-1}
\delta\left(\frac{z_{1}-z}{z_{0}}\right)\right)
a(z_{1})b(z)-{\partial \over \partial z_{0}}\left(
z_{0}^{-1}\delta\left(\frac{-z+z_{1}}{z_{0}}\right)\right)b(z)a(z_{1})
\right)\nonumber\\
&=&-{\rm Res}_{z_{1}}
\left({\partial\over \partial z_{1}} z_{0}^{-1}\delta\left(
\frac{z_{1}-z}{z_{0}}\right)\right)a(z_{1})b(z)\nonumber\\
& &+{\rm Res}_{z_{1}}\left({\partial\over \partial z_{1}}
z_{0}^{-1}
\delta\left(\frac{-z+z_{1}}{z_{0}}\right)\right)b(z)a(z_{1})\nonumber\\
&=&{\rm Res}_{z_{1}}
\left( z_{0}^{-1}\delta\left(\frac{z_{1}-z}{z_{0}}\right)a'(z_{1})
b(z)-z_{0}^{-1}
\delta\left(\frac{-z+z_{1}}{z_{0}}\right)b(z)a'(z_{1})\right)
\nonumber\\
&=&Y(a'(z),z_{0})b(z)
\end{eqnarray}
and
\begin{eqnarray}
& &[D,Y(a(z),z_{0})]b(z)\nonumber \\
&=&D(Y(a(z),z_{0})b(z))-Y(a(z),z_{0})Db(z)\nonumber \\
&=&{\partial \over \partial
z}(Y(a(z),z_{0})b(z))-Y(a(z),z_{0})b'(z)\nonumber \\
&=&{\rm Res}_{z_{1}}\left(\left({\partial \over \partial z}z_{0}^{-1}
\delta\left(\frac{z_{1}-z}{z_{0}}\right)\right)a(z_{1})b(z)
-\left({\partial \over \partial z}z_{0}^{-1}
\delta\left(\frac{z-z_{1}}{-z_{0}}\right)\right)b(z)a(z_{1})\right)
\nonumber\\
&=&{\rm Res}_{z_{1}}\left({\partial \over \partial z_{0}}\left(z_{0}^{-1}
\delta\left(\frac{z_{1}-z}{z_{0}}\right)\right)a(z_{1})b(z)
-{\partial \over \partial z_{0}}\left(z_{0}^{-1}
\delta\left(\frac{z-z_{1}}{-z_{0}}\right)\right)b(z)a(z_{1})\right)
\nonumber\\
&=&{\partial \over \partial z_{0}}Y(a(z),z_{0})b(z)\nonumber\\
&=&Y(a'(z),z_{0})b(z)\nonumber\\
&=&Y(D\cdot a(z),z_{0})b(z).\;\;\;\;\Box
\end{eqnarray}

{\bf Lemma 3.8}.  {\it Let $(M,d)$ be an object of $C^{o}$ and
let $a(z),b(z)\in F(M,d)$. Then $a(z)_{n} b(z)\in F(M,d)$.
Furthermore, if $M$ is a restricted Vir-module with
central charge $\ell$ and $a(z),b(z)$ are weak vertex operators on
$(M,L(-1))$ of weights $\alpha, \beta$, respectively, then for any
integer $n$, $a(z)_{n} b(z)$ is a weak vertex
operator of weight $(\alpha +\beta -n-1)$ on $(M,L(-1))$}.

{\bf Proof.} It is equivalent to prove the following:
\begin{eqnarray}
& &[L(-1),Y(a(z),z_{0})b(z)]={\partial \over \partial
z}(Y(a(z),z_{0})b(z));\\
& &[L(0),Y(a(z),z_{0})b(z)]\nonumber\\
&=&(\alpha +\beta)Y(a(z),z_{0})b(z)+z_{0}{\partial\over \partial
z_{0}}(Y(a(z),z_{0})b(z))+z{\partial\over \partial
z}(Y(a(z),z_{0})b(z)).
\end{eqnarray}
By definition, we have:
\begin{eqnarray}
& &{\partial\over \partial z}(Y(a(z),z_{0})b(z))\nonumber\\
&=&{\partial\over \partial z}{\rm Res}_{z_{1}}
\left(z_{0}^{-1}\delta\left(\frac{z_{1}-z}{z_{0}}\right)a(z_{1})b(z)
-z_{0}^{-1}\delta\left(\frac{z-z_{1}}{-z_{0}}\right)
b(z)a(z_{1})\right)\nonumber\\
&=&{\rm Res}_{z_{1}}\left(\left({\partial\over\partial z}z_{0}^{-1}
\delta \left(\frac{z_{1}-z}{z_{0}}\right)\right)a(z_{1})b(z)
-\left({\partial\over\partial z}z_{0}^{-1}
\delta \left(\frac{z-z_{1}}{-z_{0}}\right)\right)
b(z)a(z_{1})\right)\nonumber\\
& &+{\rm Res}_{z_{1}}
\left(\delta\left(\frac{z_{1}-z}{z_{0}}\right)a(z_{1})b'(z)
-z_{0}^{-1}\delta\left(\frac{z-z_{1}}{-z_{0}}\right)
b'(z)a(z_{1})\right)\nonumber\\
&=&-{\rm Res}_{z_{1}}
\left(\left({\partial\over \partial z_{1}}z_{0}^{-1}\delta\left(\frac
{z_{1}-z}{z_{0}}\right)\right)a(z_{1})b(z)-
\left({\partial\over \partial z_{1}}z_{0}^{-1}\delta\left(\frac{z_{1}-z}
{z_{0}}\right)\right)b(z)a(z_{1})\right)\nonumber\\
& &+{\rm Res}_{z_{1}}\left(z_{0}^{-1}\delta\left(\frac{z_{1}-z}
{z_{0}}\right)a(z_{1})b'(z)
-z_{0}^{-1}\delta \left(\frac{z_{1}-z}{z_{0}}\right)
b'(z)a(z_{1})\right)\nonumber\\
&=&{\rm Res}_{z_{1}}
\left(z_{0}^{-1}\delta\left(\frac{z_{1}-z}{z_{0}}\right)a(z_{1})b'(z)
-z_{0}^{-1}\delta \left(\frac{z-z_{1}}{-z_{0}}\right)
b'(z)a(z_{1})\right)\nonumber\\
& &-{\rm Res}_{z_{1}}
\left(z_{0}^{-1}\delta\left(\frac{z_{1}-z}{z_{0}}\right)a'(z_{1})b(z)
-z_{0}^{-1}\delta \left(\frac{z-z_{1}}{-z_{0}}\right)
b(z)a'(z_{1})\right)\nonumber\\
&=&[L(-1),Y(a(z),z_{0})b(z)] \nonumber\end{eqnarray}
and
\begin{eqnarray}
& &[L(0),Y(a(z),z_{0})b(z)]\nonumber\\
&=&{\rm Res}_{z_{1}}
\left(z_{0}^{-1}\delta\left(\frac{z_{1}-z}
{z_{0}}\right)[L(0),a(z_{1})b(z)]
-z_{0}^{-1}\delta \left(\frac{z-z_{1}}{-z_{0}}\right)
[L(0),b(z)a(z_{1})]\right)\nonumber\\
&=&{\rm Res}_{z_{1}}z_{0}^{-1}\delta\left(\frac{z_{1}-z}
{z_{0}}\right)\left(a(z_{1})[L(0),b(z)]+[L(0),a(z_{1})]b(z)\right)
\nonumber\\
& & -{\rm Res}_{z_{1}}z_{0}^{-1}\delta\left(
\frac{z-z_{1}}{-z_{0}}\right)
\left(b(z)[L(0),a(z_{1})]+[L(0),b(z)]a(z_{1})\right)\nonumber\\
&=&{\rm Res}_{z_{1}}z_{0}^{-1}\delta\left(\frac{z_{1}-z}
{z_{0}}\right)\left(\beta a(z_{1})b(z)+za(z_{1})b'(z)+\alpha
a(z_{1})b(z)+z_{1}a'(z_{1})b(z)\right)\nonumber\\
& &-{\rm Res}_{z_{1}}z_{0}^{-1}\delta\left(\frac{z-z_{1}}
{-z_{0}}\right)\left(\alpha b(z)a(z_{1})+z_{1}b(z)a'(z_{1})
+\beta b(z)a(z_{1})+zb'(z)a(z_{1})\right)\nonumber\\
&=&(\alpha +\beta )Y(a(z),z_{0})b(z)\nonumber\\
& &+{\rm Res}_{z_{1}}\left(z_{0}^{-1}\delta\left(\frac{z_{1}-z}
{z_{0}}\right)za(z_{1})b'(z)-z_{0}^{-1}\delta\left(
\frac{z-z_{1}}{-z_{0}}\right)zb'(z)a(z_{1})\right)\nonumber\\
& &-{\rm Res}_{z_{1}}\left(\left({\partial\over\partial z_{1}}z_{1}z_{0}^{-1}
\delta\left(\frac{z_{1}-z}{z_{0}}\right)\right)a(z_{1})b(z)
-\left({\partial\over\partial z_{1}}z_{1}z_{0}^{-1}
\delta\left(\frac{z-z_{1}}{-z_{0}}\right)\right)b(z)a(z_{1})\right)
\nonumber\\
&=&(\alpha +\beta )Y(a(z),z_{0})b(z)\nonumber\\
& &+{\rm Res}_{z_{1}}\left(z_{0}^{-1}\delta\left(\frac{z_{1}-z}
{z_{0}}\right)za(z_{1})b'(z)-z_{0}^{-1}\delta\left(
\frac{z-z_{1}}{-z_{0}}\right)zb'(z)a(z_{1})\right)\nonumber\\
& &-{\rm Res}_{z_{1}}(z_{0}+z)\left({\partial\over\partial z_{1}}
z_{0}^{-1}\delta\left(\frac{z_{1}-z}{z_{0}}\right)\right)a(z_{1})b(z)
\nonumber\\
& &+ {\rm Res}_{z_{1}}(z_{0}+z)\left({\partial\over
\partial z_{1}}z_{0}^{-1}
\delta\left(\frac{z-z_{1}}{-z_{0}}\right)\right)b(z)a(z_{1})
\nonumber\\
&=&(\alpha +\beta )Y(a(z),z_{0})b(z)\nonumber\\
& &+{\rm Res}_{z_{1}}\left(z_{0}^{-1}\delta\left(\frac{z_{1}-z}
{z_{0}}\right)za(z_{1})b'(z)-z_{0}^{-1}\delta\left(
\frac{z-z_{1}}{-z_{0}}\right)zb'(z)a(z_{1})\right)\nonumber\\
& &+{\rm Res}_{z_{1}}z_{0}\left(\left({\partial\over\partial z_{0}}
z_{0}^{-1}\delta\left(\frac{z_{1}-z}{z_{0}}\right)\right)a(z_{1})b(z)
-\left({\partial\over\partial z_{1}}z_{0}^{-1}
\delta\left(\frac{z-z_{1}}{-z_{0}}\right)\right)b(z)a(z_{1})\right)
\nonumber\\
& &+{\rm Res}_{z_{1}}z\left(\left({\partial\over\partial z}
z_{0}^{-1}\delta\left(\frac{z_{1}-z}{z_{0}}\right)\right)a(z_{1})b(z)
-\left({\partial\over\partial z_{1}}z_{0}^{-1}
\delta\left(\frac{z-z_{1}}{-z_{0}}\right)\right)b(z)a(z_{1})\right)
\nonumber\\
&=&(\alpha +\beta )Y(a(z),z_{0})b(z)+z{\partial\over \partial z}
\left(Y(a(z),z_{0})b(z)\right)+z_{0}{\partial\over \partial z_{0}}\left(
Y(a(z),z_{0})b(z)\right).\nonumber\;\;\;\;\Box
\end{eqnarray}

The following definition is motivated by physicists' work for example [Go].

{\bf Definition 3.9}. Two weak
vertex operators $a(z_{1})$ and $b(z_{2})$ are
said to be {\it mutually local} if there a positive integer $n$ such that
\begin{eqnarray}
(z_{1}-z_{2})^{n}a(z_{1})b(z_{2})=(z_{1}-z_{2})^{n}
b(z_{2})a(z_{1}).
\end{eqnarray}
A weak vertex operator is called a
{\it vertex operator} if it is
local with itself, a subspace $A$ of $F(M)$ is
said to be {\it local} if any two
weak vertex operators in $A$ are
mutually local, and {\it a local system} of vertex operators on $M$ is
a maximal local  subspace of $F(M)$.

{\bf Remark 3.10}. Let $V$ be a vertex algebra and let $(M, Y_{M})$ be a
$V$-module. Then the image of $V$ under the linear map
$Y_{M}(\cdot,z)$ is a local subspace of $F(M)$.

{\bf Remark 3.11}. Let $M$ be a vector space and let $a(z)$ and
 $b(z)$ be homogeneous mutually local weak vertex operators on $M$.
Let $k$ be a positive integer satisfying (3.1). Then
$a(z)_{n}b(z)=0$ whenever $n\ge k$. Thus $Y(a(z),z_{0})b(z)$
involves only finitely many negative powers of $z_{0}$. (This
corresponds to the truncation condition in the definition of a vertex operator
algebra.)

{\bf Lemma 3.12}. {\it If $a(z_{1})$ is local with $b(z_{2})$, then
$a(z_{1})$ is local with $b'(z_{2})$}.

{\bf Proof}. Let $n$ be a positive integer such that (3.17) holds. Then
\begin{eqnarray}
(z_{1}-z_{2})^{n+1}a(z_{1})b(z_{2})=(z_{1}-z_{2})^{n+1}
b(z_{2})a(z_{1}).
\end{eqnarray}
Differentiating (3.18) with respect to $z_{2}$, then using (3.17) we obtain
\begin{eqnarray}
(z_{1}-z_{2})^{n+1}a(z_{1})b'(z_{2})=(z_{1}-z_{2})^{n+1}
b'(z_{2})a(z_{1}).\;\;\;\;\Box
\end{eqnarray}

{\bf Remark 3.13}. For any vector space $M$, it follows from
Zorn's lemma that there always
exist local systems of vertex operators on $M$. Since the identity
operator $I(z)={\rm id}_{M}$ is mutually local with any weak vertex
operator on $M$, any local system contains $I(z)$.
{}From Remark 3.11 and Lemma 3.12, any local system
is closed under the derivative operator $\displaystyle{D={d\over dz}}$.

{\bf Lemma 3.14}. {\it Let $M$ be a restricted $Vir$-module with
central charge
$\ell$. Then $L(z)=\sum_{n\in {\bf Z}}L(n)z^{-n-2}$ is a (local)
vertex operator on $(M,L(-1))$ of weight two.}

{\bf Proof}. It follows from Remark 2.11 and
 Lemma 2.12 that
\begin{eqnarray}
(z_{1}-z_{2})^{k}[L(z_{1}),L(z_{2})]=0\;\;\;\mbox{ for }k\ge 4.
\end{eqnarray}
Then $L(z)$ is a local vertex operator on $M$. $\;\;\;\;\Box$

The proof of the following proposition was given by Professor Chongying Dong.

{\bf Proposition 3.15}. {\it Let $a(z)$, $b(z)$ and $c(z)$ be
 weak vertex operators on $M$. Suppose
both $a(z)$ and $b(z)$ are local with $c(z)$. Then
$a(z)_{n}b(z)$ is local with $c(z)$ for all $n\in {\bf Z}$}.

{\bf Proof}. Let $r$ be a positive integer greater than $-n$ such that the
following identities hold:
\begin{eqnarray*}
& &(z_{1}-z_{2})^{r}a(z_{1})b(z_{2})=(z_{1}-z_{2})^{r}b(z_{2})a(z_{1}),\\
& &(z_{1}-z_{2})^{r}a(z_{1})c(z_{2})=(z_{1}-z_{2})^{r}c(z_{2})a(z_{1}),\\
& &(z_{1}-z_{2})^{r}b(z_{1})c(z_{2})=(z_{1}-z_{2})^{r}c(z_{2})b(z_{1}).
\end{eqnarray*}
By definition, we have
\begin{eqnarray}
a(z)_{n}b(z)
={\rm Res}_{z_{1}}\left((z_{1}-z)^{n}a(z_{1})b(z)
-(-z+z_{1})^{n}b(z)a(z_{1})\right).
\end{eqnarray}
Since
\begin{eqnarray}
& &(z-z_{3})^{4r}\left((z_{1}-z)^{n}a(z_{1})b(z)c(z_{3})
-(-z+z_{1})^{n}b(z)a(z_{1})c(z_{3})\right)\nonumber\\
&=&\sum_{s=0}^{3r}\left(\begin{array}{c}3r\\s\end{array}\right)
(z-z_{1})^{3r-s}(z_{1}-z_{3})^{s}(z-z_{3})^{r}\cdot
\nonumber\\
& &\cdot \left((z_{1}-z)^{n}a(z_{1})b(z)c(z_{3})
-(-z+z_{1})^{n}b(z)a(z_{1})c(z_{3})\right)\nonumber\\
&=&\sum_{s=r+1}^{3r}\left(\begin{array}{c}3r\\s\end{array}\right)
(z-z_{1})^{3r-s}(z_{1}-z_{3})^{s}(z-z_{3})^{r}\cdot
\nonumber\\
& &\cdot \left((z_{1}-z)^{n}a(z_{1})b(z)c(z_{3})
-(-z+z_{1})^{n}b(z)a(z_{1})c(z_{3})\right)\nonumber\\
&=&\sum_{s=r+1}^{3r}\left(\begin{array}{c}3r\\s\end{array}\right)
(z-z_{1})^{3r-s}(z_{1}-z_{3})^{s}(z-z_{3})^{r}\cdot
\nonumber\\
& &\cdot \left((z_{1}-z)^{n}c(z_{3})a(z_{1})
b(z)
-(-z+z_{1})^{n}c(z_{3})b(z)a(z_{1})\right)\nonumber\\
&=&(z-z_{3})^{4r}\left((z_{1}-z)^{n}c(z_{3})a(z_{1})b(z)
-(-z+z_{1})^{n}c(z_{3})b(z)a(z_{1})\right),\nonumber\\
& &\mbox{}
\end{eqnarray}
we have
\begin{eqnarray}
(z-z_{3})^{4r}(a(z)_{n}b(z))c(z_{3})=
(z-z_{3})^{4r}c(z_{3})
(a(z)_{n}b(z)).\;\;\;\;\Box
\end{eqnarray}

{\bf Remark 3.16.} Let $M$ be any super vector space and let $V$ be
any local system of vertex operators on $M$. Then it follows from
Proposition 3.15, Remarks 3.11 and 3.12 and Lemmas 3.13 and 3.14
that the quadruple $(V,I(z),D, Y)$ satisfies (V1)-(V4) of Definition 2.1.

{\bf Proposition 3.17}. {\it Let $V$ be any local system of
vertex operators on $M$. Then
for any  vertex operators $a(z)$ and $b(z)$ in $V$,
$Y(a(z),z_{1})$ and $Y(b(z),z_{2})$ are mutually local
on $(V,D)$.}

{\bf Proof}. Let $c(z)$ be any weak vertex operator on
$M$. Then we have
\begin{eqnarray*}
& &Y(a(z),z_{3})Y(b(z),z_{0})c(z_{2})\\
&=&{\rm Res}_{z_{1}}
z_{3}^{-1}\delta\left(\frac{z_{1}-z_{2}}{z_{3}}\right)
a(z_{1})(Y(b(z),z_{0})c(z_{2}))\\
& &-z_{3}^{-1}\delta\left(\frac{-z_{2}+z_{1}}{z_{3}}\right)
(Y(b(z),z_{0})c(z_{2}))a(z_{1})\\
&=&{\rm Res}_{z_{1}}{\rm Res}_{z_{4}}A
\end{eqnarray*}
where
\begin{eqnarray*}
A&=& z_{3}^{-1}\delta\left(\frac{z_{1}-z_{2}}{z_{3}}\right)z_{0}^{-1}
\delta\left(\frac{z_{4}-z_{2}}{z_{0}}\right)a(z_{1})b(z_{4})c(z_{2})\\
& &-z_{3}^{-1}\delta\left(\frac{z_{1}-z_{2}}{z_{3}}\right)
z_{0}^{-1}
\delta\left(\frac{-z_{2}+z_{4}}{z_{0}}\right)a(z_{1})c(z_{2})b(z_{4})\\
& &-z_{3}^{-1}\delta\left(
\frac{-z_{2}+z_{1}}{z_{3}}\right)z_{0}^{-1}
\delta\left(\frac{z_{4}-z_{2}}{z_{0}}\right)b(z_{4})c(z_{2})a(z_{1})\\
& &+
z_{3}^{-1}\delta\left(\frac{-z_{2}+z_{1}}{z_{3}}\right)z_{0}^{-1}
\delta\left(\frac{-z_{2}+z_{4}}{z_{0}}\right)c(z_{2})b(z_{4})a(z_{1}).
\end{eqnarray*}
Similarly, we have
\begin{eqnarray}
Y(b(z),z_{0})Y(a(z),z_{3})c(z_{2})
={\rm Res}_{z_{1}}{\rm Res}_{z_{4}}B
\end{eqnarray}
where
\begin{eqnarray*}
B&=& z_{3}^{-1}\delta\left(\frac{z_{1}-z_{2}}{z_{3}}\right)z_{0}^{-1}
\delta\left(\frac{z_{4}-z_{2}}{z_{0}}\right)b(z_{4})a(z_{1})c(z_{2})\\
& &-z_{3}^{-1}\delta\left(\frac{-z_{2}+z_{1}}{z_{3}}\right)
z_{0}^{-1}
\delta\left(\frac{z_{4}-z_{2}}{z_{0}}\right)b(z_{4})c(z_{2})a(z_{1})\\
& &- z_{3}^{-1}\delta\left(\frac{z_{1}-z_{2}}{z_{3}}\right)z_{0}^{-1}
\delta\left(\frac{-z_{2}+z_{4}}{z_{0}}\right)a(z_{1})c(z_{2})b(z_{4})\\
& &+z_{3}^{-1}\delta\left(\frac{-z_{2}+z_{1}}{z_{3}}\right)z_{0}^{-1}
\delta\left(\frac{-z_{2}+z_{4}}{z_{0}}\right)c(z_{2})a(z_{1})b(z_{4}).
\end{eqnarray*}
Let $k$ be any positive integer such that
\begin{eqnarray*}
(z_{1}-z_{4})^{k}a(z_{1})b(z_{4})=(z_{1}-z_{4})^{k}b(z_{4})a(z_{1}).
\end{eqnarray*}
Since
\begin{eqnarray*}
(z_{3}-z_{0})^{k}z_{3}^{-1}\delta\left(\frac{z_{1}-z_{2}}{z_{3}}\right)z_{0}^{-1}
\delta\left(\frac{z_{4}-z_{2}}{z_{0}}\right)
=(z_{1}-z_{4})^{k}z_{3}^{-1}\delta\left(\frac{z_{1}-z_{2}}{z_{3}}\right)z_{0}^{-1}
\delta\left(\frac{z_{4}-z_{2}}{z_{0}}\right),
\end{eqnarray*}
it is clear that
locality of $a(z)$ with $b(z)$ implies the locality of $Y(a(z),z_{1})$
with $Y(b(z),z_{2}).\;\;\;\;\Box$

Now, we are ready to present our main theorem:

{\bf Theorem 3.18}. {\it Let $M$ be any
vector space and let $V$ be any local system of vertex
operators on $M$. Then $V$ is a vertex algebra and $M$ satisfies
all the conditions for module except the existence of $d$ in (M2). If
$V$ is a local system on $(M,d)$, then $(M,d)$ is  a $V$-module.}

{\bf Proof}. It follows from Proposition 2.4, Remark 3.16 and
Proposition 3.17 that
$V$ is a vertex algebra. It follows from Proposition 2.3 and
Remark 2.3 that $M$ is a
$V$-module through the linear map $Y_{M}(a(z),z_{0})=a(z_{0})$ for
$a(z)\in V. \;\;\;\;\Box$

{\bf Corollary 3.19.} {\it Let $M$ be any vector
space and let $S$ be any set of mutually local vertex operators
on $M$. Let $\langle S\rangle$ be the subspace of $F(M)$ generated by $S\cup
\{I(z)\}$
under the vertex operator multiplication (3.7) (or (3.5) for
components). Then $(\langle S\rangle , I(z), D,Y)$ is a vertex algebra with $M$
as a module.}

{\bf Proof.} It follows from Proposition 3.15 that
$\langle S\rangle$ is a local subspace of $ F(M)$. Let $A$ be a local system
containing $\langle S\rangle$ as a subspace. Then by Theorem 3.18, $A$ is a
vertex
superalgebra with $M$ as a module. Since $\langle S\rangle$ is closed under
(3.7), $\langle S\rangle$ is a vertex subalgebra. Since the ``multiplication''
(3.7) does not depend on the choice of the local system $A$,
$\langle S\rangle$ is canonical.$\;\;\;\;\Box$

{\bf Proposition 3.20}. {\it Let $M$ be a restricted $Vir$-module with
central charge $\ell$ and let $V$ be a
local system of
vertex operators on $(M, L(-1))$, containing $L(z)$. Then
the vertex operator $L(z)$
is a Virasoro element of the vertex algebra $V$.}

{\bf Proof}. First, by Theorem 3.18 $V$ is a vertex algebra with
$M$ as  a
$V$-module. Set $\omega=L(z)\in V$. By Lemma 2.7, the components of
vertex operator $Y(\omega,z_{0})$  give rise to a representation
on $V$ of central charge $\ell$ for the Virasoro algebra $Vir$.
For any
$a(z)\in V_{(h)}$, by definition we have
\begin{eqnarray}
& &L(z)_{0}a(z)=[L(-1),a(z)]=a'(z);\\
& &L(z)_{1}a(z)=[L(0),a(z)]-z^{-1}[L(-1),a(z)]=h a(z).
\end{eqnarray}
Therefore $V$ satisfies all conditions for a vertex
operator superalgebra except the requirements on the homogeneous
subspaces.$\;\;\;\;\Box$

Let $V$ be a vertex (operator) algebra and let $(M,d)$ be a
$V$-module. Then the image $\bar{V}$ of $V$ inside $F(M,d)$ is a
local subspace. By Zorn's lemma, there exists a local system
$A$ containing $\bar{V}$ as a subspace. From the vacuum property (M2)
we have:
\begin{eqnarray}
Y_{M}(\cdot,z)({\bf 1})=Y_{M}({\bf 1},z)={\rm id}_{M}=I(z).
\end{eqnarray}
For any elements $a,b\in V$, we have:
\begin{eqnarray}
& &Y_{M}(\cdot,z)(Y(a,z_{0})b)\nonumber\\
&=&Y_{M}(Y(a,z_{0})b,z)\nonumber\\
&=&{\rm Res}_{z_{1}}\left(z_{0}^{-1}\delta\left(\frac{z_{1}-z}{z_{0}}\right)
Y_{M}(a,z_{1})
Y_{M}(b,z)-z_{0}^{-1}\delta\left(\frac{z-z_{1}}
{-z_{0}}\right)Y_{M}(b,z)Y_{M}(a,z_{1})\right)\nonumber\\
&=&Y_{A}(Y_{M}(a,z),z_{0})Y_{M}(b,z).
\end{eqnarray}
Thus $Y_{M}(\cdot,z)$ is a vertex algebra homomorphism from $V$
to $A$. Conversely, let $\phi$ be a vertex algebra homomorphism from $V$
to some local system $A$ of vertex operators on $(M,d)$. Since $\phi
(u)\in A\in ({\rm End}M)[[z,z^{-1}]]$ for any $u\in V$, we use
$\phi_{z}$ for $\phi$ to indicate the dependence of $\phi(u)$ on $z$.
For any formal variable $z_{1}$, set

$\phi_{z_{1}}(a)=\phi_{z}(a)|_{z=z_{1}}$ for any $a\in V$.
We define
$Y_{M}(a,z)u=\phi_{z} (a)$ for $a\in V$. By definition we have:
\begin{eqnarray}
Y_{M}({\bf 1},z)=\phi_{z}({\bf 1})=I(z)={\rm id}_{M}.
\end{eqnarray}
For any elements $a,b\in V$, we have:
\begin{eqnarray}
& &Y_{M}(Y(a,z_{0})b,z_{2})\nonumber\\
&=&\phi_{z} (Y(a,z_{0})b)|_{z=z_{2}}\nonumber\\
&=&\left(Y_{A}(\phi_{z}(a),z_{0})\phi_{z}(b)\right)|_{z=z_{2}}\nonumber\\
&=&{\rm Res}_{z_{1}}\left(z_{0}^{-1}\delta\left(\frac{z_{1}-z_{2}}
{z_{0}}\right)\phi_{z_{1}}(a)
\phi_{z_{2}}(b)-z_{0}^{-1}\delta\left(\frac{z_{2}-z_{1}}
{-z_{0}}\right)\phi_{z_{2}}(b)\phi_{z_{1}}(a)\right)\nonumber\\
&=&{\rm Res}_{z_{1}}\left(z_{0}^{-1}\delta\left(\frac{z_{1}-z_{2}}{z_{0}}
\right)Y_{M}(a,z_{1})
Y_{M}(b,z_{2})-z_{0}^{-1}\delta\left(\frac{z_{2}-z_{1}}
{-z_{0}}\right)Y_{M}(b,z_{2})Y_{M}(a,z_{1})\right).\nonumber\\
& &\mbox{ }
\end{eqnarray}
It follows from Remark 2.4 that $(M,d,Y_{M})$ is a $V$-module.
Therefore, we have proved:

{\bf Proposition 3.21.} {\it Let $V$ be a vertex (operator)
algebra. Then
giving a $V$-module $(M,d)$ is equivalent to giving a vertex algebra
homomorphism from $V$ to some local system of vertex operators on $(M,d)$.}

\section{Vertex operator algebras and modules
associated to some infinite-dimensional Lie algebras}
In this section, we shall use the machinery we built in Section 3 to
study vertex operator algebras and modules associated
to the representations for some well-known infinite-dimensional Lie algebras
such as the Virasoro algebra and affine Lie algebras. Most of the material
 presented in this section was taken from [Li2]. (See [FZ] for a different
approach.)

Let us start with an abstract result which will be used in this
section.

{\bf Proposition 4.1}. {\it Let $(V,D, {\bf 1},Y)$ be a vertex
algebra and let $(M,d,Y_{M})$ be a $V$-module. Let
$u\in M$ such that $du=0$. Then the linear map}
\begin{eqnarray}
f: V\rightarrow M; a\mapsto a_{-1}u\;\;\;\mbox{ {\it for} }a\in V,
\end{eqnarray}
{\it is a $V$-homomorphism.}

{\bf Proof.} It follows from the proof of Proposition 3.4 [Li1].$\;\;\;\;\Box$

For any complex numbers $c$ and $h$, let $M(c,h)$ be the Verma module
for the Virasoro algebra $Vir$ with central charge $c$ and with
lowest weight $h$. Let ${\bf 1}$ be a lowest weight vector of $M(c,0)$.
Then $L(-1){\bf 1}$ is a singular vector, i.e., $L(n)L(-1){\bf 1}=0$ for
$n\ge 1$. Set $\bar{M}(c,0)=M(c,0)/\langle L(-1){\bf 1}\rangle$, where $\langle
L(-1){\bf
1}\rangle$ denotes the submodule of $M(c,0)$ generated by $L(-1){\bf 1}$.
Denote by
$L(c,h)$ the (unique) irreducible quotient module of $M(c,h)$. By
slightly abusing notations, we still use ${\bf 1}$ for the image of
${\bf 1}$ for both $\bar{M}(c,0)$ and $L(c,0)$.

{\bf Proposition 4.2}. {\it For any complex number $c$,
$\bar{M}(c,0)$ has a natural vertex operator algebra structure and any
restricted $Vir$-module $M$ of central charge $c$ is a weak
$\bar{M}(c,0)$-module. In
particular, for any complex number $h$, $M(c,h)$ is a $\bar{M}(c,0)$-module.}

{\bf Proof.}
Let $M$ be any restricted $Vir$-module with central charge $c$. Then
$\bar{M}(c,0)\oplus M$ is a restricted $Vir$-module. By Lemma 3.14,
$L(z)$ is a local vertex operator on $(\bar{M}(c,0)\oplus M,L(-1))$.
Then by Corollary 3.19, $V=\langle L(z)\rangle$ is a vertex algebra
with $\bar{M}(c,0)\oplus M$ as a module. Consequently, both
$\bar{M}(c,0)$ and $M$ are $V$-modules. By Lemma 2.5, the components
of $Y(L(z),z_{0})$ on $V$ satisfy the Virasoro relation. Since
$L(z)_{n}I(z)=0$ for $n\ge 0$,
$V$ is a lowest weight $Vir$-module with lowest weight $0$, so that $V$
is a quotient module of $\bar{M}(c,0)$. Let ${\bf 1}$ be
a lowest weight vector of $\bar{M}(c,0)$. Since
$L(-1){\bf 1}=0$, by Proposition 4.1, we have a $V$-homomorphism from $V$ to
$\bar{M}(c,0)$ mapping $I(z)$ to ${\bf 1}$. Then it follows that $V$
is isomorphic to $\bar{M}(c,0)$. Therefore $\bar{M}(c,0)$ is a
vertex operator algebra and any restricted $Vir$-module $M$ is a weak module.
$\;\;\;\;\Box$

{\bf Remark 4.3.} It follows that $L(c,0)$ is a quotient vertex
operator algebra of $\bar{M}(c,0)$.

Let $({\bf g},B)$ be a pair consisting of a finite-dimensional
Lie algebra ${\bf g}$ and a nondegenerate symmetric invariant bilinear form
$B$ on ${\bf g}$.
Set
\begin{eqnarray}
\tilde{\bf g}={\bf C}[t,t^{-1}]\otimes {\bf g}\oplus {\bf C}c.
\end{eqnarray}
Then we define
\begin{eqnarray}
[a_{m},b_{n}]=[a,b]_{m+n}+m\delta_{m+n,0}\langle a,b\rangle c,\;\;[c, x_{m}]=0
\end{eqnarray}
for any $a,b, x\in {\bf g}$ and
$m,n\in {\bf Z}$, where $x_{m}$ stands for $t^{m}\otimes x$.
For any $x\in {\bf g}$, we set
\begin{eqnarray}
x(z)=\sum_{n\in {\bf Z}}x_{n}z^{-n-1}.
\end{eqnarray}
Then the defining relations (4.3) of $\tilde{{\bf g}}$ are
equivalent to the following equations:
\begin{eqnarray}
& &[a(z_{1}),b(z_{2})]=z_{1}^{-1}\delta\left(\frac{z_{2}}{z_{1}}\right)
[a,b](z_{2})+
z_{1}^{-2}\delta'\left(\frac{z_{2}}{z_{1}}\right)\langle
a,b\rangle c,\\
& &[x(z),c]=0\;\;\;\mbox{ for any }a,b,x\in {\bf g}, m,n\in
{\bf Z}.
\end{eqnarray}
Set
\begin{eqnarray}
N_{+}=t{\bf C}[t]\otimes {\bf g},\; N_{-}=t^{-1}{\bf C}[t^{-1}]\otimes {\bf
g},\;N_{0}={\bf g}\oplus {\bf C}c.
\end{eqnarray}
Then we obtain a triangular decomposition $\tilde{{\bf g}}=N_{+}\oplus
N_{0}\oplus N_{-}$. Let $P=N_{+}+N_{0}$ be the parabolic subalgebra.
For any ${\bf g}$-module $U$ and any complex
number $\ell$, denote by $M_{({\bf g},B)}(\ell,U)$ the generalized Verma
module [Le] or Weyl
module with $c$ acting as scalar $\ell$. Namely,
$M_{({\bf g},B)}(\ell,U)=U(\tilde{{\bf g}})\otimes_{U(P)}U$.
For any $\tilde{{\bf g}}$-module $M$, we may consider $x(z)$ for $x\in
{\bf g}$ as an element of $({\rm End}M)[[z,z^{-1}]]$. Recall that a
$\tilde{{\bf g}}$-module $M$ is said to be {\it restricted} if for any
$u\in M$, $(t^{k}{\bf C}[t]\otimes {\bf g})u=0$ for $k$ sufficiently
large. Then a $\tilde{{\bf g}}$-module $M$
is restricted if and only if $x(z)$ for all $x\in {\bf g}$ are weak
vertex operators on $M$.

{\bf Theorem 4.4.} {\it For any complex number $\ell$, $M_{({\bf
g},B)}(\ell,{\bf C})$ has a natural vertex algebra structure and any
restricted $\tilde{{\bf g}}$-module $M$ of level $\ell$ is a $M_{({\bf
g},B)}(\ell,{\bf C})$-module.}

{\bf Proof.} Let $M$ be any restricted $\tilde{{\bf g}}$-module of
level $\ell$. Then $W=M_{({\bf g},B)}(\ell,{\bf C})\oplus M$ is also a
restricted  $\tilde{{\bf g}}$-module of level $\ell$.
It follows from  Lemma 2.3 and (4.5)-(4.6) that
$\bar{{\bf g}}=\{a(z)|a\in {\bf g}\}$ is a local subspace of
$F(W)$. Let $V$ be the subspace of $F(W)$ generated by
all $\bar{\bf g}$. Then by Corollary 3.19, $V$ is a vertex
superalgebra and $W$ is a $V$-module. Consequently, both $M$ and
$M_{({\bf g},B)}(\ell,{\bf C})$ are $V$-modules.
It follows from Lemma 2.5 and
(4.5)-(4.6) that $V$ is a $\tilde{{\bf g}}$-module (of
level $\ell$) with a vector $I(z)$ satisfying $P\cdot I(z)=0$, so that $V$ is a
quotient $\tilde{{\bf g}}$-module of $M_{({\bf g},B)}(\ell,{\bf C})$.

To finish the proof, we only need to prove that $V$ is isomorphic to
$M_{({\bf g},B)}(\ell,{\bf C})$ as a $V$-module. Let $d$ be the
endomorphism of $M_{({\bf g},B)}(\ell,{\bf C})$ such
that
\begin{eqnarray}
d\cdot {\bf 1}=0,\;[d, a_{m}]=-ma_{m-1}\;\;\mbox{ for }a\in {\bf g}.
\end{eqnarray}
Then $[d,a(z)]=a'(z)$ for any $a\in {\bf g}$. Then $(M_{({\bf
g},B)}(\ell,{\bf C}),d)$ is a $V$-module. It
follows from Proposition 4.1 and the universal property of $M_{({\bf
g},B)}(\ell,{\bf C})$ that $V$ and $M_{({\bf g},B)}(\ell,{\bf C})$ are
isomorphic $V$-modules. $\;\;\;\;\Box$

{\bf Remark 4.6.} It is clear that $M_{({\bf
g},B)}(\ell,{\bf C})=M_{({\bf g},\alpha B)}(\alpha^{-1}\ell,{\bf C})$ for any
nonzero complex number $\alpha$ (see for example [Lia]).

\end{document}